\documentstyle[amssymb]{amsart}

\numberwithin{equation}{section}

\theoremstyle{plain}
	\newtheorem{thm}{Theorem}[section]
	\newtheorem{lem}[thm]{Lemma}
	
	\newtheorem{prop}[thm]{Proposition}
\theoremstyle{definition}

\theoremstyle{remark}
	\newtheorem{rem}[thm]{Remark}

\newcommand{\NC}[3]{C_{{#1},{#2}}({#3};q)}
\newcommand{\comment}[1]{}
\title{Raising operators of row type \\ for Macdonald polynomials }
\author{Yasushi KAJIHARA and Masatoshi NOUMI}
\address{Department of Mathematics, Kobe University, 
Rokko, Kobe 657-8501, Japan}
\email{kaji@@math.s.kobe-u.ac.jp, noumi@@math.s.kobe-u.ac.jp}

\begin{document}

\maketitle
\begin{center}
{\small
Department of Mathematics, Kobe University
}
\end{center}
\begin{abstract}
We construct certain raising operators of row type
for Macdonald's symmetric polynomials by an interpolation method.
\end{abstract}

\section{Introduction}

Throughout this paper, we denote by $J_\lambda(x;q,t)$ the integral form of 
Macdonald's symmetric polynomial in $n$ variables $x=(x_1,\ldots,x_n)$ 
(of type $A_{n-1}$) associated with a partition $\lambda$ (\cite{macdo1}). 
For each $m=0,1,2,\ldots$, we consider a $q$-difference operator $B_m$ 
which should satisfy the following condition: 
{\sl For any partition $\lambda=(\lambda_1,\lambda_2,\ldots)$ whose longest 
part $\lambda_1$ has length $\le m$, one has}
\begin{equation} \label{rr}
B_m J_\lambda (x;q,t) =
\begin{cases}
 J_{(m,\lambda)}(x;q,t) & \text{if} \quad \ell(\lambda)<n, \\ 
0 & \text{if} \quad \ell(\lambda)=n,
\end{cases}
\end{equation}
{\sl where $(m,\lambda)=(m,\lambda_1,\lambda_2,\ldots)$ stands for the 
partition obtained by adding a row of length $m$ to $\lambda$. }
An operator $B_m$ having this property will be called a 
{\em raising operator of row type\/}
for Macdonald polynomials. 
With such operators, 
the Macdonald polynomial $J_\lambda(x;q,t)$ 
for a general partition $ \lambda = (\lambda_1,
\lambda_2 ,\ldots , \lambda_n ) $
can be expressed as 
\begin{equation}
B_{\lambda_1} B_{\lambda_2} \ldots B_{\lambda_n}.1 = J_\lambda (x;q,t)
\quad (\lambda_1 \ge  \lambda_2 \ge \ldots \ge \lambda_n \ge 0).
\end{equation}
Namely, one can obtain $ J_\lambda ( x ; q ,t ) $ by an 
successive application of the operators $ B_m $ 
starting from  $J_\phi(x;q,t)=1 $. 

The purpose of this paper is to give an explicit construction 
of such operators $B_m$ $(m=0,1,2,\ldots)$. 
These operators $B_m$ can be considered as a {\em dual version} of 
the raising operators {\em of column type} 
introduced by A.N.~Kirillov and the second author \cite{kn1}, \cite{kn2}. 
We remark that, as to the Hall-Littlewood polynomials 
(the case when $q=0$), such a class of raising operators 
$B_m$ of row type has been implicitly employed in Macdonald \cite{macdo1}, 
Chapter III, (2.14): 
\begin{equation}
B_m=(1-t)\sum_{i=1}^n x_i^m 
\left(\prod_{j\ne i}\frac{x_i-tx_j}{x_i-x_j}\right)\,T_{0,x_i}
\end{equation}
for $m=1,2,\ldots$, where $T_{0,x_i}$ is the ``0-shift operator '' in $x_i$, 
namely, the substitution of zero for $x_i$.  
Our raising operators of row type for Macdonald polynomials can be 
considered as a generalization of these operators for 
Hall-Littlewood polynomials.

\medskip
We will propose first a theorem of unique existence for raising operators
of row type.
For each multi-index $\alpha=(\alpha_1,\ldots,\alpha_n)\in {\Bbb N}^n$,
we set $|\alpha|=\alpha_1+\cdots+\alpha_n$ 
and 
\begin{equation}
x^\alpha=x_1^{\alpha_1}\cdots x_n^{\alpha_n}, \quad
T_{q,x}^\alpha=T_{q,x_1}^{\alpha_1}\cdots T_{q,x_n}^{\alpha_n}, 
\end{equation}
where $T_{q,x_i}$ is the $q$-shift operator in $x_i$, defined by 
\begin{equation}
T_{q,x_i}f\,(x_1,\ldots,x_i,\ldots,x_n)= f(x_1,\ldots,qx_i,\ldots,x_n)
\end{equation}
for $i=1,\ldots,n$.

\begin{thm}\label{ut}
For each $ m=0,1,2,\ldots$, there exists a unique $q$-difference operator 
\begin{equation}
B_m = \sum_{|\gamma| \le m} b^{(m)}_\gamma (x) T_{q,x}^\gamma  
\end{equation}
of order $ \le m $ satisfying the condition \eqref{rr}, where 
$b^{(m)}_\gamma(x)$ are rational functions in $x$ with coefficients 
in ${{\Bbb Q}}(q,t)$. 
Furthermore, the operator $B_m$ is invariant under the action of 
the symmetric group ${\frak S}_n$ of degree $n$. 
\end{thm}

We will also determine the operator $B_m$ explicitly by an  
interpolation method.  
In the following, we use the notation $\alpha\le \beta$ for the 
partial ordering of multi-indices defined by
\begin{equation}\label{po-mi}
\alpha\le\beta\quad \Longleftrightarrow\quad 
\alpha_i\le \beta_i\quad (i=1,\ldots,n). 
\end{equation}
In order to describe the coefficients of our raising operators, 
we introduce a variant of {\em $q$-binomial coefficients} 
$\NC{\alpha}{\beta}{x}$ including the variables $x=(x_1,\ldots,x_n)$. 
For any pair $(\alpha,\beta)$ of multi-indices such that $\alpha\ge\beta$,
we set
\begin{align}\label{n-coef}
\NC{\alpha}{\beta}{x}&=\prod_{1\le i,j\le n}
\frac{(q^{\alpha_i-\beta_j+1}x_i/x_j)_{\beta_j}}
{(q^{\beta_i-\beta_j+1}x_i/x_j)_{\beta_j}}
\\
&=\prod_{j=1}^n
\frac{(q^{\alpha_j-\beta_j+1})_{\beta_j}}
{(q)_{\beta_j}}
\prod_{i\ne j}
\frac{(q^{\alpha_i-\beta_j+1}x_i/x_j)_{\beta_j}}
{(q^{\beta_i-\beta_j+1}x_i/x_j)_{\beta_j}}
\notag
\end{align}
with the notation 
$(a)_k=(a;q)_k=(1-a)(1-a q)\cdots(1-aq^{k-1})$ of the $q$-shifted 
factorial.  
We remark that, if $n=1$, $\NC{\alpha}{\beta}{x}$ reduce to the ordinary 
$q$-binomial coefficients 
$\left[\begin{matrix} \alpha \\ \beta \end{matrix}\right]_q$. 
\begin{thm}\label{xf}
The $q$-difference operator $B_m$ of Theorem \ref{ut} can be expressed 
in the form
\begin{equation}
B_m =  \sum_{|\alpha| = m } b_\alpha^{(m)} (x) \, \phi_\alpha^{(m)}
 (x;T_{q,x}), 
\end{equation}
where
\begin{align}
b_\alpha^{(m)}(x) 
&=
(-1)^{|\alpha|}
q^{\sum_i { \alpha_i  \choose 2}} x^\alpha \sum_{\beta \le \alpha }
(-1)^{|\beta |} q^{ |\beta| \choose 2} 
\NC{\alpha}{\beta}{x} 
\\
&\qquad\cdot
\prod_{i,j=1}^n \frac{
(tq^{-\beta_j+1}x_i/x_j)_{\beta_j}
(q^{-\alpha_j+1}x_i/x_j)_{\alpha_j-\beta_j}
}{
(q^{\alpha_i-\alpha_j+1}x_i/x_j)_{\alpha_j}
}
\notag
\end{align}
and
\begin{equation}
\phi^{(m)}_{\alpha}(x;T_{q,x})=
\sum_{\beta\le \alpha}
(-1)^{|\alpha| - |\beta|}
 q^{|\alpha|- |\beta| + 1 \choose 2 } 
\NC{\alpha}{\beta}{x}\, T_{q,x}^\beta
\end{equation}
for each $\alpha$ with $|\alpha|=m$. 
\end{thm}

In the course of the proof of Theorem \ref{xf}, we will make use of 
a variant of the {\em $q$-binomial theorem} for our $\NC{\alpha}{\beta}{x}$, 
which might also deserve attention 
(see Proposition \ref{prop-nqbinom} in Section 5).
\begin{thm}
For any $\alpha\in {\Bbb N}^n$, one has 
\begin{equation}\label{qbf}
\sum_{\beta\le\alpha}(-1)^{|\beta|}q^{|\beta| \choose 2} 
\NC{\alpha}{\beta}{x} u^{|\beta|}=(u)_{|\alpha|}. 
\end{equation}
\end{thm}
\noindent
We remark that formula (\ref{qbf}) also implies a generalization of 
{\em $q$-Chu-Vandermonde formulas} 
\begin{equation}
\sum_{\beta\le\alpha, |\beta|=r} 
\prod_{j=1}^n \left[
\begin{matrix}\alpha_j\\ \beta_j \end{matrix}
\right]_q
\prod_{i\ne j}
\frac{(q^{\alpha_i-\beta_j+1}x_i/x_j)_{\beta_j}}
{(q^{\beta_i-\beta_j+1}x_i/x_j)_{\beta_j}}
=
\left[
\begin{matrix}n \\ r \end{matrix}
\right]_q
\end{equation}
for any $\alpha$ with $|\alpha|=n$ and $0\le r\le n$. 
 
\medskip
After recalling some basic facts about Macdonald polynomials in Section 2, 
we will prove the uniqueness and the existence of raising operators of 
row type in Section 3 and in Section 4, respectively.  
Explicit formulas for the $q$-difference operators 
$\phi^{(m)}_\alpha(x;T_{q,x})$ and the coefficients $b^{(m)}_{\alpha}(x)$ 
($|\alpha|=m$) of Theorem \ref{xf} 
will be given in Section 5 and in Section 6, respectively. 


\section{Macdonald Polynomials}

In order to fix the notation, 
we recall some basic facts about Macdonald's 
symmetric polynomials of type  $A_{n-1}$. For the details
see \cite{macdo1}.

\medskip

Let ${\Bbb K}[x]={\Bbb K}[x_1,x_2,\ldots,x_n]$ be the ring of polynomials 
in $n$ variables $x=(x_1, x_2, \ldots , x_n)$ with coefficients 
in ${\Bbb K}={\Bbb Q}(q,t)$, and 
${\Bbb K}[x]^{{\frak S}_n}$ the subring of all invariant polynomials under
the natural action of the symmetric group ${\frak S}_n$ of 
degree $n$. 

Macdonald's commuting family of $q$-difference operators 
$D_1,D_2,\ldots,D_n$ is defined by the generating function
\begin{eqnarray}
D_x (u ; q, t ) &=& \sum_{r = 0 }^n ( - u )^r D_r \\
&=& \sum_{ K \subset \{ 1, \ldots , n\} } ( - u )^{|K|} 
q^{ |K| \choose 2 }\prod_{ i \in K, j \notin K} 
\frac{1 - t x_i / x_j }{1 - x_i / x_j } \prod_{ i \in K }  
T_{ q, x_i}. \nonumber
\end{eqnarray}
Note that $D_x(u;q,t)$ has the determinantal formula
\begin{eqnarray}
D_x (u ; q, t ) &=&
\frac{1}{ \Delta(x) }
\det ( x_j^{ n-i } ( 1 - u t^{ n-i }T_{ q, x_i}) )_{i,j} \\
&=& 
\frac{1}{ \Delta(x) } \sum_{ w \in \frak S_n}
\epsilon ( w ) w ( \prod_{i = 1 }^n x_i^{n-i}  
( 1 - u t^{ n-i }T_{ q, x_i} )), \nonumber
\end{eqnarray}
where $\Delta(x) = \prod_{i < j} ( x_i - x_j ) $ .
Macdonald's symmetric polynomials 
$ P_\lambda (x)=P_\lambda(x;q,t)$ are the joint eigenfunctions 
of the operators $D_1,\ldots,D_n$ on 
${\Bbb K}[x]^{{\frak S}_n}$, satisfying the equations 
\begin{equation}
D_x (u ) P_\lambda (x) = P_\lambda (x) 
\prod_{i=1 }^n ( 1 - u q^{ \lambda_i} t^{n-i});
\end{equation}
each $P_\lambda(x)$ is 
normalized so that the coefficient of $x^\lambda$ should be equal to $1$.
The integral form $ J_\lambda (x) = J_\lambda (x ; q,t)$ 
of $P_\lambda(x)$ 
is defined as 
\begin{equation}
J_\lambda (x ; q,t) = c_\lambda P_\lambda (x ; q,t) ,\quad 
c_\lambda = \prod_{s \in \lambda } ( 1 - q^{a(s)} t^{l(s)+1 }). 
\end{equation}
It is known in fact that $J_\lambda(x)$ are linear combinations of 
monomial symmetric functions with 
coefficients in ${\Bbb Z}[q,t]$ (see \cite{kn1} for example). 
\par
We recall 
that the Macdonald polynomials have the generating function
\begin{equation}\label{gf}
\prod_{i=1 }^n \prod_{j=1 }^m \,  ( 1 + x_i y_j ) = \sum_{\lambda } 
P_\lambda (x;q,t)  P_{\lambda'} (y;t,q), 
\end{equation}
for another set of variables $y=(y_1,\ldots,y_m)$, 
where $\lambda' $ stands for the conjugate partition of $\lambda $,
and the summation is taken over all partitions $\lambda$ such that 
$l( \lambda ') = \lambda_1 \le m$, $l(\lambda ) = {\lambda'}_1\le n  $.
This formula will be the key to our study 
of raising operators of row type. 
Notice that the dual version of the generation function \eqref{gf} 
has been employed in \cite{kn1} for the construction of 
raising operators of column type. 


\section{Raising operators of row type and their uniqueness}

%
Fixing a nonnegative integer $m$, 
we will prove in this section the uniqueness of a 
$q$-difference operator 
\begin{equation}\label{dop}
B_m = \sum_{|\gamma| \le m} b^{(m)}_\gamma (x) T_{q,x}^\gamma
\quad(b^{(m)}_\gamma(x)\in {\Bbb K}(x))
\end{equation}
of order $\le m$ such that 
\begin{equation}\label{rop}
B_m J_\lambda (x;q,t) = 
	\begin{cases}
	J_{ (m,\lambda ) } 
    (x;q,t) \quad& \mbox{if}\quad l (\lambda ' ) \le m , l(\lambda) < n, \\
	0 &  \mbox{if}\quad l(\lambda ' ) \le m , l(\lambda) = n,
	\end{cases}
\end{equation}
where $ (m,\lambda ) = ( m, \lambda_1, \lambda_2, \ldots)$.
We remark that the invariance of 
$B_m$ under the action of ${\frak S}_n $ follows immediately 
from the uniqueness theorem.  
Existence of such an operator will be established in the next section. 

\medskip

\begin{lem}\label{lem-key-id}
A $q$-difference operator $ B_m $ of order $\le m$ 
in the form \eqref{dop} satisfies the 
condition \eqref{rop} if and only if the following equality holds: 
\begin{equation}\label{rop-g}
   B_{m , x } \prod_{ i=1}^n \prod_{ j=1}^m ( 1 + x_i y_j )
= \frac{1}{ y_1 \ldots y_m } D_y ( 1 ; t, q) \prod_{i = 1}^n 
\prod_{ j=1}^m ( 1 + x_i y_j ).
\end{equation}
\end{lem}

\begin{pf}
Note first that, for each  partition $ \mu=(\mu_1,\ldots,\mu_m) $ 
of length $\le m$, one has 
\begin{equation}\label{psi-pol}
 \frac{1}{ y_1 \ldots y_m } D_y ( 1 ; t, q)  P_{\mu }( y ; t, q ) = 
	\begin{cases} 
             	P_{\mu - (1)^m }( y ; t, q)
		\prod_{i = 1 }^m ( 1 - q^{ m- i} t^{ \mu_i } ) 
		& \text{if $\mu_m>0$, } \\
		 0 &  \text{if $ \mu_m = 0 $.}
	\end{cases}
\end{equation}
Hence we obtain
\begin{eqnarray}\label{mac-g}
&&
\frac{1}{ y_1 \ldots y_m } D_y ( 1 ; t, q) \prod_{i = 1}^n \prod_{j= 1}^m
( 1 + x_i y_j )  \\    
&=&
\sum_{ l(\nu ) \le n , l(\nu') =m }
P_{\nu} ( x ; q,t ) 
P_{\nu'-(1)^{m} } ( y ; t,q )
\prod_{i = 1 }^m ( 1 - q^{ m- i} t^{ (\nu ' )_i} )\nonumber
\\
&=& \sum_{ l(\lambda ) \le n-1 , l'(\lambda ) \le m }
P_{(m, \lambda )} ( x ; q,t ) 
P_{\lambda ' } ( y ; t,q )
\prod_{i = 1 }^m ( 1 - q^{ m- i} t^{ (\lambda ' )_i  + 1} )\nonumber.
\end{eqnarray}
This implies that 
equation (\ref{rop-g}) is equivalent to the condition
\begin{eqnarray}
B_m P_\lambda ( x ; q,t)  
&=& 
\begin{cases}
	0 & \text{( if $l(\lambda  ) = n  $ )}  \\
	P_\lambda ( x ; q,t) \prod_{i=1}^m(1-q^{m-i}t^{(\lambda')_i+1})
	& \text{( if $l(\lambda  ) < n $ )}
\end{cases}
\end{eqnarray}
for any $\lambda$ with $l(\lambda')\le m$. 
It is easily seen that this coincides with condition (\ref{rop}) in terms of 
the integral forms. 
\end{pf} 

By making the action of $D_y(1;t,q)$ in (\ref{rop-g}) explicit, we obtain 

\begin{prop}
A $q$-difference operator $B_m$ of order $\le m$ is a raising operator 
of row type for Macdonald polynomials if and only if its coefficients 
satisfy the following identity of rational functions:
\begin{eqnarray}\label{key-id}
&&
\sum_{|\gamma|\le m } b^{(m)}_\gamma (x) 
\prod_{i = 1}^n \prod_{j= 1}^m \frac{1 +  q^{\gamma_i }x_i y_j}{1+x_iy_j} \\
&=&
\frac{1}{ y_1 \ldots y_m }
\sum_{ K \subset \{ 1, \ldots , n \}}
( - 1)^{|K|} q^{|K| \choose 2 }
\prod_{ k \in K , l \notin K }
\frac{ 1 - q y_k /y_l }{ 1 - y_k /y_l }
\prod_{i = 1}^n \prod_{k \in K }\frac{ 1 + t x_i y_k }{1+x_i y_k} \nonumber.
\end{eqnarray}
\end{prop}

\begin{rem}
By the determinantal representation of $ D_y (1 ; t, q ) $, 
equality (\ref{key-id}) can also be rewritten in the form 
\begin{eqnarray}
&&
\sum_{|\gamma|\le m } b^{(m)}_\gamma (x) 
\prod_{i = 1}^n \prod_{j= 1}^m
\frac{1 +  q^{\gamma_i }x_i y_j}{1+x_i y_j}\\
&=&\frac{1}{ y_1 \ldots y_m  \Delta (y)}  
\det \left( y_j^{m-i}\big( 1 - q^{m-i} \prod_{r =1}^n
\frac{ 1 + t x_r y_j }{ 1 + x_r y_j }\big)
\right)_{i,j}. \nonumber
\end{eqnarray}
\end{rem}
 
Let now $B$ and $B'$ be two $q$-difference operators of order $\le m$ and 
suppose that they both satisfy the condition (\ref{rop}) of raising operators. 
Then by Lemma \ref{lem-key-id} one has 
\begin{equation}
(B_x-B'_x)\prod_{i=1}^n \prod_{j=1}^m (1+x_i y_j)=0. 
\end{equation}
Hence the uniqueness of  $B_m$ of Theorem \ref{ut} follows immediately 
from the following general proposition on $q$-difference operators.

\begin{prop}\label{qdop}
Let $P=\sum_{|\gamma|\le m} a_\gamma(x) T_{q,x}^\gamma $ be a $q$-difference 
operator of order $\le m$ with coefficients in ${\Bbb K}(x)$.  
\newline 
$(a)$ If $P_x \prod_{i=1}^n \prod_{j=1}^m (1+x_i y_j)=0$, then $P=0$ as a
$q$-difference  operator. 
\newline
$(b)$ If $Pf(x)=0$ for any symmetric polynomial 
$f(x)\in{\Bbb K}[x]^{\frak S_n}$
of  degree $\le mn$, then $P=0$ as a $q$-difference operator. 
\end{prop}
\noindent
Since the statement (b) follows from (a), 
we give a proof of (a) of Proposition.
For each multi-index $\alpha\in{\Bbb N}^n$ 
with $|\alpha|=m$, we define a point 
$p_\alpha(x)\in {\Bbb K(x)}^m$ by
\begin{eqnarray}\label{def-p}
p_\alpha(x)&=&(-1/x_1, -1/qx_1,\ldots,-1/q^{\alpha_1-1}x_1,\ldots,
\\
& &\quad -1/x_n, -1/qx_n,\ldots,-1/q^{\alpha_{n}-1}x_n). \nonumber
\end{eqnarray}
Then we have 
\begin{lem}\label{daiji}
For any multi-index $\gamma\in {\Bbb N}^n$, one has 
\begin{eqnarray}
\prod_{i=1}^n \prod_{j=1}^m (1+q^{\gamma_i}x_i y_j)\big|_{y=p_\alpha(x)}
&=&
\prod_{i=1}^n \prod_{j=1}^n
\prod_{\nu=0}^{\alpha_j-1}(1-q^{\gamma_i-\nu}x_i/x_j)\\
&=&
\prod_{1\le i,j\le n}(q^{\gamma_i-\alpha_j+1}x_i/x_j)_{\alpha_j}.\nonumber
\end{eqnarray}
In particular, one has 
$\prod_{i=1}^n \prod_{j=1}^m (1+q^{\gamma_i}x_i y_j)\big|_{y=p_\alpha(x)}=0$ 
unless $\gamma\ge \alpha$.
\end{lem}
Under the assumption of Proposition \ref{qdop},(a), we may assume that 
$a_\alpha(x)\ne 0$ for some $\alpha\in{\Bbb N}^n$ 
with $|\alpha|=m$ without loosing 
generality.  (If $P$ is of order $l<m$, 
set $y_{l+1}=\ldots=y_m=0$ and apply the following argument 
by replacing $m$ by $l$.) 
The assumption on $P$ implies
\begin{equation}
\sum_{|\gamma|\le m} a_\gamma(x) 
\prod_{i=1}^n \prod_{j=1}^m (1+q^{\gamma_i}x_i y_j)=0.
\end{equation}
Evaluating this equality at $y=p_\alpha(x)$, we have
\begin{equation}
a_\alpha(x) \prod_{1\le i,j\le n}(q^{\alpha_i-\alpha_j+1}x_i/x_j)_{\alpha_j}=0
\end{equation}
by Lemma \ref{daiji},
since, if $|\gamma|\le m$ and $\gamma\ge \alpha$, then $\gamma=\alpha$. 
This contradicts to the assumption $a_\alpha(x)\ne0$. 
This completes the proofs of Proposition \ref{qdop} and the uniqueness 
of $B_m$ in Theorem \ref{ut}.

\section{Existence of $B_m$}

In this section, we discuss the existence of a raising operator $B_m$. 

\medskip

We begin with a lemma which will play an important role in the 
following argument. 
\begin{lem}\label{intp}
Let $ F(y) \in {\Bbb K} (x) [y]^{{\frak S}_m} $ be a symmetric polynomial 
in $y=(y_1,\ldots,y_m)$ with coefficients in ${\Bbb K}(x)$, and suppose 
that $F(y)$ is of degree $ \le  n-1 $ in $y_j$ for each $j=1,\ldots,m$. 
If $ F(p_\alpha(x))=0 $ for all $\alpha\in{\Bbb N}^n$ with $|\alpha|=m$, 
then $F(y)$ is identically zero as a polynomial in $y$. 
\end{lem}
\begin{pf}
We prove Lemma by the induction on $m$. 
The case when $m = 1 $ is obvious since $F(y)$ is of degree $\le n-1$ and 
has $n$ distinct zeros $-1/x_1,\ldots,-1/x_n$. 
For $m\ge 2$, we first expand $ F(y)$ in terms of $ y_m$ as follows: 
\begin{equation}
 F( y ) = F( y_1 , \ldots , y_m ) = \sum_{ i=0 }^{n-1}  
F_i( y_1 , \ldots , y_{m - 1 } ) y_m^i,
\end{equation}
where each coefficient $F_i ( y_1 , \ldots , y_{m - 1 } ) $ has 
degree $\le n-1$ 
in all $y_j$ ($j=1,\ldots,m-1$). 
Let $\beta\in{\Bbb N}^n$ a multi-index with $|\beta|=m-1$
and consider the polynomial 
\begin{equation}
 f( y_m ) = F(p_\beta(x),y_m ) = \sum_{ i=0 }^{n-1}  
F_i( p_\beta(x)) y_m^i,
\end{equation}
by evaluating $F(y)$ at $(y_1,\ldots,y_{m-1})=p_\beta(x)$. 
From the assumption on $ F(y)$, it follows that the polynomial 
$f(y_m)$  has $n$ 
distinct zeros $ y_m = 
- 1/ q^{\beta_i} x_i $ ($i=1,\ldots,n$).
Hence $ f( y_m )  $ is identically $ 0 $ as a polynomial in $y_m$.
This implies that  
 $F_i( p_\beta(x))=0 $ for each $ i = 0 ,\ldots , m-1$ and 
for any $\beta\in{\Bbb N}^n$ 
with $|\beta|=m-1$. 
By the induction hypothesis, we conclude 
that the coefficients $ F_i( y_1 , \ldots , y_{m - 1 } )$ are 
identically zero 
as polynomials in $(y_1,\ldots,y_{m-1})$, namely, 
$F(y)$ is identically zero as a polynomial in $y=(y_1,\ldots,y_m)$. 
\end{pf}

In view of Lemma \ref{lem-key-id}, we propose to construct a 
$q$-difference operator 
\begin{equation}\label{bexp}
B=\sum_{|\alpha|\le m} b_\alpha(x)T_{q,x}^\alpha
\end{equation}
of order $\le m$ such that 
\begin{equation} \label{key-id2}
   B_{x } \prod_{ i=1}^n \prod_{ j=1}^m ( 1 + x_i y_j )
= \frac{1}{ y_1 \ldots y_m } D_y ( 1 ; t, q) \prod_{i = 1}^n 
\prod_{ j=1}^m ( 1+x_i y_j ).
\end{equation} 
In the following, we denote the left-hand side and the right-hand side 
of this 
equality by $\Phi(x;y)$ and by $\Psi(x;y)$, respectively. 
In terms of the coefficients $b_\alpha(x)$, $\Phi(x;y)$ is expressed as  
\begin{equation}\label{phib}
\Phi(x;y)=\sum_{|\alpha|\le m}b_\alpha(x)
\prod_{i=1}^n\prod_{j=1}^m(1+q^{\alpha_i}x_i y_j).
\end{equation}
Note also that $\Psi(x;y)$ is a polynomial in $y=(y_1,\ldots,y_m)$ and 
has degree $\le n-1$ in each $y_j$ ($j=1,\ldots,m$) as can be 
seen from (\ref{psi-pol}).  
Hence, by Lemma \ref{intp}, we see that 
 $B$ satisfies the desired equality if and only if 
\begin{enumerate}
\item $\Phi(x;y)$ is of degree $\le n-1$ in each $y_j$ for $j=1,\ldots,m$.
\item $\Phi(x;p_\alpha(x))=\Psi(x;p_\alpha(x))$ for all $\alpha\in{\Bbb N}^n$ 
with $|\alpha|=m$. 
\end{enumerate}


\medskip
Suppose now that the operator $B$  has the property (1) mentioned above. 
Since the degree of $\Phi(x;y)$ in $y_j$ is less than $n$ 
for each $j=1,\ldots,m$,  we have 
\begin{equation} 
\Phi ( x; y ) \prod_{i=1}^n  \prod_{j=1}^m ( 1 + x_i y_j )^{ - 1 }
|_{ y_1 \to \infty , \ldots ,  y_m \to \infty } = 0. 
\label{y,phi}
\end{equation}
Hence by (\ref{phib}) we obtain 
\begin{equation}
\sum_{|\alpha | \le m} b_\alpha  (x) q ^{ |\alpha | m } = 0, \quad
\text{i.e.,}\quad
b_0(x)=-\sum_{0<|\alpha | \le m} b_\alpha  (x) q ^{ |\alpha | m }.
\end{equation}
This implies that $ B $ can be represent as
\begin{equation} \label{b1}
B =  \sum_{1\le |\alpha | \le m } b_\alpha (x) 
(T_{q,x}^\alpha  - q ^{ |\alpha| m } ).
\end{equation}
Note that a general $B$ of order $\le m$ has an expression of this form
if and only if
\begin{equation}
F_1(x;y_1)=\Phi (x;y) \prod_{i=1}^n  \prod_{j=1}^m (1+x_i y_j)^{-1}
|_{y_2\to\infty,\ldots, y_m\to\infty} 
\end{equation}
is of degree $\le n-1$ in $y_1$. 
We now show inductively that, 
for $l =0,1,\ldots,m$, $B$ can be represented as follows: 
\begin{eqnarray}\label{bla}
B &=& \sum_{l   \le |\alpha| \le m} b_\alpha  (x) \phi_{l;\alpha}
 (x, T_{q,x}),
\end{eqnarray}
where 
\begin{equation}\label{blb}
   \phi_{l;\alpha} (x, T_{q,x}) = T_{q,x}^\alpha 
   + \sum_{\beta <  \alpha , \, |\beta | < l}  \phi_{l;\alpha , \beta }(x)
 T_{q,x}^\beta.
\end{equation}
Assume that we have  constructed such an expression for $l$ with $l<m$.
Note that
\begin{eqnarray} 
\Phi ( x; y ) &=&
\sum_{l \le |\alpha | \le m} b_\alpha  (x) \Bigg(
\prod_{i=1}^n  \prod_{j=1}^m ( 1 +  q^{\alpha_i } x_i y_j ) \\
&&\quad 
+ \sum_{ \beta \le \alpha,\, |\beta | < l} 
\phi_{l;\alpha, \beta }(x)
\prod_{i=1}^n  \prod_{j=1}^m ( 1 +  q^{\beta_i } x_i y_j )\Bigg). \nonumber
\end{eqnarray} 
Since property (1) of $\Phi(x;y)$ implies 
\begin{equation} 
\Phi ( x; y ) \prod_{i=1}^n  \prod_{j=l+1}^m ( 1 + x_i y_j )^{ - 1 }
|_{ y_{l +1} \to \infty , \ldots ,  y_m \to \infty } = 0,
\end{equation}
we obtain the relation
\begin{eqnarray}\label{phi-l}  
&&\sum_{l \le |\alpha | \le m} b_\alpha  (x) \Bigg(q^{|\alpha | ( m -l ) }
\prod_{i=1}^n  \prod_{j=1}^l ( 1 +  q^{\alpha_i } x_i y_j ) 
\\
&&
+ \sum_{ \beta \le \alpha, \, |\beta | < l}  
\phi_{l;\alpha, \beta }(x) q^{|\beta | ( m -l ) }
\prod_{i=1}^n  \prod_{j=1}^l ( 1 +  q^{\beta_i } x_i y_j ) \Bigg) = 0. 
\nonumber
\end{eqnarray} 
In this formula 
we consider to specialize $y'=(y_1,\ldots,y_l)$ at $p_\gamma(x)$,
with the notation of  (\ref{def-p}), for each $\gamma$ with $|\gamma|=l$.
By Lemma \ref{daiji}, 
$\displaystyle{\prod_{i=1}^n 
\prod_{j=1}^l( 1 + q^{\beta_i} x_i y_j )}|_{y'=p_\gamma(x)} =0$ 
unless $\beta \ge  \gamma$.  
Hence formula (\ref{phi-l}) with $y'=p_\gamma(x)$ gives rise to 
\begin{eqnarray} 
&&
 b_\gamma (x) q^{l ( m -l ) }
\prod_{1\le i,j\le n} (q^{\gamma_i -\gamma_j+1} x_i/ x_j )_{\gamma_j}
\\
&&
+ \sum_{|\alpha | > l} b_\alpha  (x) q^{|\alpha | ( m -l ) }
\prod_{1\le i,j\le n} (q^{\alpha_i -\gamma_j+1} x_i/ x_j )_{\gamma_j}= 0 
\nonumber.
\end{eqnarray}
From this we have 
\begin{equation}
b_\gamma (x) = - \sum_{ \alpha > \gamma} b_\alpha(x)  
\psi_{\alpha, \gamma}(x),
\end{equation}
where 
\begin{eqnarray}\label{psi}
\psi_{\alpha, \gamma}(x)&=&
q^{(|\alpha|-|\gamma|)(m-|\gamma|)}\prod_{1\le i,j\le n}
\frac{(q^{\alpha_i-\gamma_j+1}x_i/x_j)_{\gamma_j}}
{(q^{\gamma_i-\gamma_j+1}x_i/x_j)_{\gamma_j}}\\
&=&q^{(|\alpha|-|\gamma|)(m-|\gamma|)} \NC{\alpha}{\gamma}{x}\nonumber 
\end{eqnarray}
with the notation of (\ref{n-coef}).
Note that $  \psi_{\alpha, \gamma}(x)$ depends on $m$ 
but does {\em not} on $B$.
Thus we obtain
\begin{eqnarray}
B &=& \sum_{|\gamma | = l } b_\gamma (x) \phi_{l;\gamma}( x, T_{q,x} )
+\sum_{l < |\alpha | \le m } b_\alpha (x) 
\phi_{l;\alpha }( x, T_{q,x} ) \\
&=& \sum_{l+1 \le  |\alpha | \le m } b_\alpha (x) 
\phi_{l+1;\alpha}( x, T_{q,x} ) \nonumber.
\end{eqnarray}
where $ \phi_{l+1;\alpha}( x, T_{q,x} ) $ 
($l+1\le |\alpha|\le m$) are determined by
\begin{equation}\label{rec-phi0}
\phi_{l+1;\alpha}( x, T_{q,x} )
= \phi_{l;\alpha}( x, T_{q,x} )-\sum_{\gamma <\alpha , |\gamma| = l } 
\psi_{\alpha,\gamma}(x)\phi_{l;\gamma}( x,T_{q,x} ).
\end{equation}
In other words, the coefficients of $\phi_{l+1;\alpha}(x;T_{q,x})$ 
are determined by the recurrence formula 
\begin{equation}\label{rec-phi}
\phi_{l+1;\alpha,\beta}(x)=
\phi_{l;\alpha,\beta}(x)-
\sum_{\beta< \gamma <\alpha,\, |\gamma|=l}
\psi_{\alpha,\gamma}(x) \phi_{l;\gamma,\beta}(x)
\end{equation}
for all $\beta$ such that $\beta<\alpha$ and $|\beta|<l$. 
In this induction procedure, it is also seen by Lemma \ref{intp}
that a general $B$ of order $\le m$ has an expression of this form 
(\ref{bla}) with (\ref{blb}) if and only if
\begin{equation}
F_l(x;y_1,\ldots,y_m)=
\Phi ( x; y ) \prod_{i=1}^n  \prod_{j=1}^m ( 1 + x_i y_j )^{ - 1 }
|_{ y_{l+1} \to \infty , \ldots ,  y_m \to \infty } 
\end{equation}
is of degree $\le n-1$ in $y_j$ for each $j=1,\ldots,l$. 

\medskip
In this way,  we can define the $q$-difference operators 
$\phi_{l;\alpha}(x;T_{q,x})$ ($l\le|\alpha|\le m$) 
for $l=0,\ldots,m$, inductively 
on $l$ by (\ref{rec-phi0}).  
Note that these operators depend on the $m$ that we have fixed in advance, 
but do {\em not} on the operator $B$. 
By using the operators we obtained at the final step $l=m$, 
we have the expression 
\begin{equation}
B=\sum_{|\alpha|=m} b_\alpha(x) \phi^{(m)}_\alpha(x;T_{q,x})
\end{equation}
for $B$, where $\phi^{(m)}_\alpha(x;T_{q,x})=\phi_{m;\alpha}(x;T_{q,x})$ .

\medskip
From this construction, 
we obtain the following proposition.
\begin{prop}\label{binphi}
For  each $ \alpha \in {\Bbb N}^n$ with $ |\alpha | = m $, define the 
$q$-difference operator $\phi^{(m)}_{\alpha}(x;T_{q,x})$ as above. 
Then, for any $q$-difference operator $B$ of order $\le m$ 
with coefficients in ${\Bbb K}(x)$, 
the following two conditions are equivalent.
\begin{quote}
$(a)$ $\Phi(x;y)=B_x \prod_{i=1}^n  \prod_{j=1}^m ( 1 + x_i y_j ) $ is of 
degree $\le n-1$  in $ y_j$ for each $j=1,\ldots,m$.
\\
$(b)$ $B$ is represented as 
\begin{equation}
B = \sum_{|\alpha| = m } b_\alpha (x) \phi_\alpha^{(m)} (x, T_{q,x})
\end{equation}
for some $b_\alpha (x) \in {\Bbb K} (x)$.
\end{quote}
\end{prop}

We now consider a $q$-difference operator $B$ of the form 
Proposition \ref{binphi}, (b), 
so that $\Phi(x;y)$ is of degree $\le n-1$ in each $y_j$ ($j=1,\ldots,m$). 
With $\Psi(x;y)$ being the right-hand side of (\ref{key-id2}),
the equality $\Phi(x;y)=\Psi(x;y)$ holds if and only if 
$\Phi(x;p_\alpha(x))=\Psi(x;p_\alpha(x))$ for any $\alpha$ with $|\alpha|=m$,
as we remarked before. 
Since
\begin{equation}
\Phi(x;p_\alpha(x))=
b_\alpha(x) \prod_{1\le i,j\le n}(q^{\alpha_i-\alpha_j+1}x_i/x_j)_{\alpha_j}
\end{equation}
by Lemma \ref{daiji}, the coefficients $b_\alpha(x)$ are determined as
\begin{equation}\label{b-by-psi}
b_\alpha(x)=\Psi(x;p_\alpha(x))
\prod_{1\le i,j\le n}(q^{\alpha_i-\alpha_j+1}x_i/x_j)_{\alpha_j}^{-1}
\end{equation}
for all $\alpha$ with $|\alpha|=m$.
This completes the proof of existence of a raising operator $B_m$. 

\medskip
From the recurrence formula (\ref{rec-phi})  we see that,
for any $\alpha$ with $l\le|\alpha|\le m$,  
the coefficients $\phi_{l;\alpha,\beta}(x)$ of $\phi_{l;\alpha}(x;T_{q,x})$ 
are expressed as 
\begin{equation}
\phi_{l;\alpha,\beta}(x)
=\sum_{r=1}^l (-1)^r 
\sum_{\alpha>\gamma_1>\ldots>\gamma_r=\beta;\, |\gamma_1|<l}
\psi_{\alpha,\gamma_1}(x)\psi_{\gamma_1,\gamma_2}(x)
\cdots \psi_{\gamma_{r-1},\gamma_r}(x)
\end{equation}
for all $\beta$ with $\beta<\alpha, |\beta|<l$.
In particular, we have 
\begin{prop} \label{prop-phi}
For any pair $(\alpha,\beta)$ of multi-indices with $\beta\le\alpha$, define 
a rational function $\psi^{(m)}_{\alpha,\beta}(x)$ by 
\begin{eqnarray}
\psi^{(m)}_{\alpha,\beta}(x)&=&
q^{(|\alpha|-|\beta|)(m-|\beta|)}\,\NC{\alpha}{\beta}{x}
\\
&=&q^{(|\alpha|-|\beta|)(m-|\beta|)}\prod_{1\le i,j\le n}
\frac{(q^{\alpha_i-\beta_j+1}x_i/x_j)_{\beta_j}}
{(q^{\beta_i-\beta_j+1}x_i/x_j)_{\beta_j}}.\nonumber
\end{eqnarray}
Then, 
for any $\alpha\in {\Bbb N}^n$ with $|\alpha|=m$, 
the coefficients of the $q$-difference operator 
\begin{equation}
\phi^{(m)}_{\alpha}(x;T_{q,x})=\sum_{\beta\le\alpha}
\phi^{(m)}_{\alpha,\beta}(x)T_{q,x}^\beta
\end{equation}
are determined by the formula
\begin{equation}
\phi^{(m)}_{\alpha,\beta}(x)
=\sum_{r=0}^m (-1)^r
\sum_{\alpha=\gamma_0>\gamma_1>\ldots>\gamma_r=\beta}
\psi^{(m)}_{\gamma_0,\gamma_1}(x)
\cdots \psi^{(m)}_{\gamma_{r-1},\gamma_r}(x),
\end{equation}
where the summation is taken over all paths in the lattice ${\Bbb N}^n$ 
connecting $\alpha$ and $\beta$. 
\end{prop}
In the next section, we will give explicit formulas for these coefficients
$\phi^{(m)}_{\alpha,\beta}(x)$.


\section{Explicit formulas for $\phi^{(m)}_\alpha(x;T_{q,x})$}

The goal of this section is to give the explicit formula
\begin{equation}\label{xf-phi1}
\phi^{(m)}_{\alpha}(x;T_{q,x})=\sum_{\beta\le\alpha}
(-1)^{|\alpha|-|\beta|} q^{|\alpha|-|\beta|+1 \choose 2}
\, \NC{\alpha}{\beta}{x} \, T_{q,x}^\beta
\end{equation}
for 
$\phi^{(m)}_\alpha (x,T_{q,x} )$ $(|\alpha|=m)$ as in Theorem \ref{xf}. 
With the notation of Proposition \ref{prop-phi}, 
this formula is equivalent to
\begin{eqnarray}\label{xf-phi2}
\phi^{(m)}_{\alpha,\beta}(x)&=&(-1)^{|\alpha|-|\beta|} 
q^{|\alpha|-|\beta|+1 \choose 2}
\, \NC{\alpha}{\beta}{x}\\
&=&(-1)^{|\alpha|-|\beta|} q^{|\alpha|-|\beta|+1 \choose 2}
\prod_{1\le i,j\le n}
\frac{(q^{\alpha_i-\beta_j+1}x_i/x_j)_{\beta_j}}
{(q^{\beta_i-\beta_j+1}x_i/x_j)_{\beta_j}}.\nonumber
\end{eqnarray}
for $\beta\le\alpha$.

\medskip
In view of the dependence of
$\psi^{(m)}_{\alpha,\beta}(x)$ on $m$ (see Proposition \ref{prop-phi}), 
we define a function $g_{\alpha,\beta}(x)$ by
\begin{equation}\label{def-g}
g_{\alpha,\beta}(x)=q^{-(|\alpha|-|\beta|)|\beta|}\,\NC{\alpha}{\beta}{x}
\end{equation}
for any $\alpha, \beta \in{\Bbb N}^n$ with $\beta\le\alpha$, 
so that 
$\psi^{(m)}_{\alpha,\beta}(x)=q^{(|\alpha|-|\beta|)m}g_{\alpha,\beta}(x)$.
With these $g_{\alpha,\beta}(x)$, 
we also define a function $f_{\alpha,\beta}(x)$
by
\begin{equation}\label{def-f}
f_{\alpha,\beta}(x)=\sum_{r=0}^{|\alpha|-|\beta|} (-1)^r
\sum_{\alpha=\gamma_0>\gamma_1>\ldots>\gamma_r=\beta}
g_{\gamma_0,\gamma_1}(x)
\cdots g_{\gamma_{r-1},\gamma_r}(x)
\end{equation}
for any $\alpha, \beta \in{\Bbb N}^n$ with $\beta\le\alpha$. 
Then by Proposition \ref{prop-phi} we have
\begin{equation}
\phi^{(m)}_{\alpha,\beta}(x)=q^{(|\alpha|-|\beta|)m}
f_{\alpha,\beta}(x)
\end{equation}
if $|\alpha|=m$ and $\beta\le\alpha$. 
Hence, the formula (\ref{xf-phi2}) follows from the following proposition. 

\begin{prop}\label{xf-f}
Define the rational functions $f_{\alpha,\beta}(x)$ ($\beta\le\alpha$)
by the formulas 
$(\ref{def-f})$  together with  $(\ref{def-g})$.
Then they can be determined as 
\begin{equation}
f_{\alpha,\beta}(x)=(-1)^{|\alpha|-|\beta|} 
q^{-{|\alpha|-|\beta|\choose 2}-(|\alpha|-|\beta|)|\beta|}\,
\NC{\alpha}{\beta}{x}
\end{equation}
for any $\alpha,\beta$ with $\beta\le\alpha$.
\end{prop}

For the proof of Proposition \ref{xf-f}, notice that 
the functions $f_{\alpha,\beta}(x)$ are defined as the matrix elements 
of the inverse matrix of 
the lower unitriangular matrix $G=(g_{\alpha,\beta}(x))_{\alpha,\beta}$. 
Hence we have only to show the inverse matrix of $G$ is given by 
$G^{-1}=(\widetilde{f}_{\alpha,\beta}(x))_{\alpha,\beta}$ with 
\begin{equation}
\widetilde{f}_{\alpha,\beta}(x)=(-1)^{|\alpha|-|\beta|} 
q^{-{|\alpha|-|\beta|\choose 2}-(|\alpha|-|\beta|)|\beta|}\,
\NC{\alpha}{\beta}{x}. 
\end{equation}
Proposition \ref{xf-f} thus reduces to 
\begin{lem}\label{inv-g}
For any $\alpha,\beta$ with $\alpha>\beta$, one has 
\begin{equation}
\sum_{\alpha\ge \gamma\ge \beta} 
\widetilde{f}_{\alpha,\gamma}(x) \, g_{\gamma,\beta}(x)
=0. 
\end{equation}
\end{lem}
By the definition of $g_{\alpha,\beta}(x)$ and 
$\widetilde{f}_{\alpha,\beta}(x)$,
we have
\begin{eqnarray} \label{mat^(-1)}
&&
\sum_{\alpha \ge \gamma \ge \beta }  \widetilde{f_{\alpha , \gamma }} (x)
g_{\gamma,\beta } (x)\\
&=&
\sum_{\alpha \ge \gamma \ge \beta }
(-1)^{|\alpha| - |\gamma | }
q^{-{|\alpha| - |\gamma | \choose 2 }- ( |\alpha| - |\gamma |)  |\gamma |
-( |\gamma| - |\beta |)  |\beta | } \,
\NC{\alpha}{\gamma}{x}\,\NC{\gamma}{\beta}{x}.\nonumber
\end{eqnarray}
Just as in the case of binomial coefficients, 
it is directly shown that our $\NC{\alpha}{\beta}{x}$
satisfy the following identity:
\begin{eqnarray}
\NC{\alpha}{\gamma}{x}\,\NC{\gamma}{\beta}{x}
&=&\NC{\alpha}{\beta}{x}\,
\prod_{i,j}
\frac{ (q^{\gamma_i - \beta_j + 1} x_i/x_j)_{\alpha_i - \gamma_i}}
{ (q^{\gamma_i - \gamma_j + 1} x_i/x_j)_{\alpha_i - \gamma_i}}
\\
&=&\NC{\alpha}{\beta}{x}\,
\NC{\alpha-\beta}{\alpha-\gamma}{1/q^\alpha x}
\nonumber
\end{eqnarray}
where $1/q^\alpha x=(1/q^{\alpha_1}x_1,\ldots,1/q^{\alpha_n}x_n)$. 
Hence we obtain
\begin{eqnarray} 
&&
\sum_{\alpha \ge \gamma \ge \beta }  \widetilde{f_{\alpha , \gamma }}(x)
g_{\gamma,\beta } (x)
=
q^{-(|\alpha|-|\beta|)|\beta|} \,\NC{\alpha}{\beta}{x}
\\
&&\quad\cdot
\sum_{\alpha \ge \gamma \ge \beta }
(-1)^{|\alpha| - |\gamma | }
q^{-{|\alpha|-|\gamma|\choose 2}- (|\alpha|-|\gamma|)(|\gamma|-|\beta|)}
\NC{\alpha-\beta}{\alpha-\gamma}{1/q^\alpha x}.\nonumber
\end{eqnarray}
Setting $ \alpha - \beta = \lambda$
and $\alpha - \gamma = \mu $,
the last summation can be rewritten in the form 
\begin{equation}\label{gamma}
\sum_{ 0\le\mu\le\lambda}
( -1 )^{|\mu |} q^{| \mu |( 1 - |\lambda |)} q^{ | \mu | \choose  2 } 
C_{\lambda,\mu}(1/q^\alpha x).
\end{equation}
Hence Lemma \ref{inv-g} is reduced to proving that this formula becomes zero.
It is in fact a special case of 
the following analogue of the $q$-binomial theorem. 
(Replace $x$ by $1/q^\alpha x$ and set $u=q^{1-|\lambda|}$ in 
(\ref{n-bi1}) below, to see that (\ref{gamma}) becomes zero.) 
\begin{prop}\label{prop-nqbinom}
For any $ \lambda \in{\Bbb N}^n$, one has 
\begin{equation} \label{n-bi1}
\sum_{0 \le \mu \le \lambda}
( -u )^{|\mu |}q^{ | \mu | \choose  2 } 
\NC{\lambda}{\mu}{x}=( u )_{|\lambda |},
\end{equation}
where $u$ is an indeterminate.
\end{prop}
\begin{pf}
This ``$q$-binomial theorem'' follows from an identity for Macdonald's 
$q$-difference operator $D_z(u;t,q)$ in $N$ variables 
$z=(z_1,\ldots,z_N)$ with  $N=|\lambda|$.
Since $D_z (u;t,q). 1 =  (u)_N$, we have
\begin{equation}\label{mac0}
\sum_{K\subset \{1,\ldots,N\}} (-u)^{|K|}q^{ |K|\choose 2} 
\prod_{k\in K; l\notin K} \frac{1-q z_k/z_l}{1-z_k/z_l} = (u)_N.
\end{equation}
For a multi-index $\lambda\in {\Bbb N}^n$ with $|\lambda|=N$,
let us specialize (\ref{mac0}) at $ z=p_\lambda(x)$
with the notation of (\ref{def-p}). 
Note that, when we specialize $z$ at $p_\lambda(x)$, the indexing set 
$\{1,\ldots,N\}$ is divided into $n$ blocks with cardinality 
$\lambda_1,\ldots,\lambda_n$, respectively.
Furthermore, for a configuration $K$ of points in $\{1,\ldots,N\}$,
the product $ \prod_{k\in K; l\notin K}(1-q z_k/z_l)/(1-z_k/z_l)$ 
becomes zero unless the elements of $K$ should be packed to the left 
in each block.  
Such configurations $K$ are parameterized by multi-indices 
$\mu\le\lambda$ such that $|\mu|=|K|$ and 
that $\mu_i$ denotes the number of points of $K$ 
sitting in the $i$-th block for $i=1,\ldots,n$.  
For such a $K$, one has
\begin{eqnarray}
\prod_{k\in K; l\notin K} 
\frac{1-q z_k/z_l}{1-z_k/z_l} \bigg|_{z=p_\lambda(x)}
&=&\prod_{1\le i,j\le n} \prod_{\mu_i\le a<\lambda_i; 0\le b<\mu_j}
\frac{1-q^{a-b+1}x_i/x_j}{1-q^{a-b}x_i/x_j}\\
&=&\prod_{1\le i,j\le n}
\frac{(q^{\lambda_i-\mu_j+1}x_i/x_j)_{\mu_j}}
{(q^{\mu_i-\mu_j+1}x_i/x_j)_{\mu_j}} 
=\NC{\lambda}{\mu}{x}. \nonumber
\end{eqnarray}
(The indices are renamed by $k\to(j,b)$, $l\to(i,a)$.)
Hence we obtain (\ref{n-bi1}).
\end{pf}
This completes the proof of formula (\ref{xf-phi1}).

\begin{rem}
In the case of one variable, 
equation (\ref{n-bi1}) reduces the ordinary $q$-binomial theorem
\begin{equation}
\sum_{k=0 }^l (-1)^k q^{ k \choose 2} u^k  
\left[\begin{matrix} l \\ k \end{matrix}\right]_q  
= ( u )_l.
\end{equation}
If we take the coefficient of $u^k$ in formula (\ref{n-bi1}), we obtain 
\begin{equation}
\sum_{ \mu \le \lambda , | \mu| = k }
\prod_{j = 1}^n  \left[\begin{matrix}
\lambda_j \\ \mu_j
\end{matrix} \right]_q 
\prod_{i \neq j}\frac
{ (q^{\lambda_i  -\mu_j   +1} x_i/ x_j )_{\mu_j}}
{ (q^{\mu_i  -\mu_j   +1} x_i/x_j )_{\mu_j}} =
\left[\begin{matrix}
|\lambda| \\ k 
\end{matrix}\right]_q,
\end{equation}
for $k=0,1,\ldots,|\lambda|$.
This gives a generalization of the $q$-Chu-Vandermonde formula.
From (\ref{n-bi1}), we also obtain another type of 
$q$-Chu-Vandermonde formula for our $\NC{\alpha}{\beta}{x}$:
\begin{equation}
\sum\Sb
\mu\le \alpha, \nu\le \beta \\ 
|\mu |+ |\nu| = k \endSb
q^{(|\alpha| + |\mu|)|\nu|}
\NC{\alpha}{\mu}{x} \NC{\beta}{\nu}{x}
=
\left[\begin{matrix}
|\alpha| + |\beta| \\  k
\end{matrix}\right]_q.
\end{equation}
\end{rem}

\section{Determination  of $ b^{(m)}_\alpha(x) $ }

We have already proved that our raising operator 
\begin{equation}
B_m=\sum_{|\gamma|\le m} b^{(m)}_\gamma(x) T_{q,x}^\gamma
\end{equation}
of row type for Macdonald polynomials has an expression
\begin{equation}
B_m=\sum_{|\alpha|= m} b^{(m)}_\alpha(x) \phi^{(m)}_\alpha(x;T_{q,x}),
\end{equation}
with the $q$-difference operators $\phi^{(m)}_\alpha(x;T_{q,x})$ of 
(\ref{xf-phi1}).
In this section, we give explicit formulas 
for $ b^{(m)}_\alpha(x) $ for all $\alpha $ 
with $ | \alpha | = m $.

\par\medskip
As we already remarked in Section 4, 
the coefficients $b^{(m)}_\alpha(x)$ ($|\alpha|=m$) 
are determined by 
\begin{equation}\label{b-by-psi2}
b_\alpha(x)=\Psi(x;p_\alpha(x))
\prod_{1\le i,j\le n}(q^{\alpha_i-\alpha_j+1}x_i/x_j)_{\alpha_j}^{-1},
\end{equation}
where 
\begin{equation}
\Psi(x;y)=\frac{1}{y_1\ldots y_n} D_y(1;t,q)\prod_{i=1}^n 
\prod_{j=1}^m(1+x_i y_j).
\end{equation}
(See (\ref{b-by-psi}).)
Recall that 
\begin{eqnarray}
\Psi(x;y)
&=&
\frac{1}{y_1\cdots y_m} \sum_{K\in \{ 1, \ldots m\}} 
(-1)^{|K|} q^{|K| \choose 2}
\prod_{k \in K, l \notin K} \frac{1- q y_k/y_l}{1- y_k/y_l}\nonumber  \\
&&
\prod_{i=1}^n \bigg\{
\prod_{k \in K} ( 1 + t x_i y_k)
\prod_{ l \notin K} ( 1 + x_i y_l) \bigg\}. \nonumber 
\end{eqnarray}
We specialize this formula at $y=p_\alpha(x)$ 
for each $\alpha$ with $|\alpha|=m$, 
in the same way as we did in the proof of Proposition \ref{prop-nqbinom}.  
All the subsets $K$ that give rise to nonzero summands 
after the specialization 
$y=p_\alpha(x)$ are parameterized by the multi-indices $\beta$ such that 
$\beta\le\alpha$ and $|\beta|=K$.  
With this parameterization, we already showed that 
\begin{equation}
\prod_{k \in K, l \notin K} 
\frac{1- q y_k/y_l}{1- y_k/y_l}\bigg|_{y=p_\alpha(x)}
=\NC{\alpha}{\beta}{x}.
\end{equation}
Renaming the indices by $k\to(j,b)$, we have 
\begin{eqnarray}
&&\prod_{i=1}^n \bigg\{
\prod_{k \in K} ( 1 + t x_i y_k)
\prod_{ l \notin K} ( 1 + x_i y_l) \bigg\} \\
&&=
\prod_{1\le i,j\le n} \prod_{b=0}^{\beta_j-1}(1-tq^{-b}x_i/x_j)
\prod_{b=\beta_j}^{\alpha_j-1}(1-q^{-b}x_i/x_j) \nonumber\\
&&=
\prod_{1\le i,j\le n} (tq^{-\beta_j+1}x_i/x_j)_{\beta_j}
(q^{-\alpha_j+1}x_i/x_j)_{\alpha_j-\beta_j}.
\nonumber
\end{eqnarray}
Hence we have 
\begin{eqnarray}
\Psi(x;p_\alpha(x))
&=&(-1)^m q^{\sum_{i}{\alpha_i \choose 2}} x^\alpha 
\sum_{\beta\le \alpha} (-1)^{|\beta|} q^{|\beta| \choose 2}
\NC{\alpha}{\beta}{x}
\nonumber  \\
&&
\prod_{1\le i,j\le n} (tq^{-\beta_j+1}x_i/x_j)_{\beta_j}
(q^{-\alpha_j+1}x_i/x_j)_{\alpha_j-\beta_j}.
\nonumber
\end{eqnarray}
By (\ref{b-by-psi2}), we finally obtain
\begin{eqnarray}
b^{(m)}_\alpha(x)
&=&q^{\sum_{i}{\alpha_i \choose 2}} x^\alpha 
\sum_{\beta\le \alpha} (-1)^{|\alpha|-|\beta|} q^{|\beta| \choose 2}
\NC{\alpha}{\beta}{x}
\nonumber  \\
&&\qquad\cdot
\prod_{1\le i,j\le n} 
\frac{(tq^{-\beta_j+1}x_i/x_j)_{\beta_j}
(q^{-\alpha_j+1}x_i/x_j)_{\alpha_j-\beta_j}}
{(q^{\alpha_i-\alpha_j+1}x_i/x_j)_{\alpha_j}}
\nonumber \\
&=&  q^{\sum_{i}{\alpha_i \choose 2}} x^\alpha 
\sum_{\beta\le \alpha} (-1)^{|\alpha|-|\beta|} 
q^{{|\beta| \choose 2}}
\nonumber\\
&&\quad\cdot
\prod_{1\le i,j\le n} 
\frac{(tq^{-\beta_j+1}x_i/x_j)_{\beta_j}
(q^{-\alpha_j+1}x_i/x_j)_{\alpha_j-\beta_j}}
{(q^{\beta_i-\beta_j+1}x_i/x_j)_{\beta_j}
(q^{\alpha_i-\alpha_j+1}x_i/x_j)_{\alpha_j-\beta_j}},
\nonumber
\end{eqnarray}
for any $\alpha$ with $|\alpha|=m$.
This completes the proof of Theorem \ref{xf}.


\end{document}